\documentclass{article}%
\usepackage{amsfonts}
\usepackage{graphicx}
\usepackage{amsmath}
\usepackage{amssymb}%
\providecommand{\U}[1]{\protect\rule{.1in}{.1in}}
\providecommand{\U}[1]{\protect\rule{.1in}{.1in}}
\providecommand{\U}[1]{\protect\rule{.1in}{.1in}}
\newtheorem{theorem}{Theorem}

\newtheorem{corollary}[theorem]{Corollary}

\newtheorem{definition}[theorem]{Definition}

\newtheorem{lemma}[theorem]{Lemma}
\newtheorem{notation}[theorem]{Notation}

\newtheorem{proposition}[theorem]{Proposition}
\newtheorem{remark}[theorem]{Remark}

\newenvironment{proof}[1][Proof]{\textbf{#1.} }{\ \rule{0.5em}{0.5em}}
\begin{document}

\title{Differential Geometry of Microlinear Fr\"{o}licher Spaces II}
\author{Hirokazu Nishimura\\Institute of Mathematics, University of Tsukuba\\Tsukuba, Ibaraki, 305-8571, Japan}
\maketitle

\begin{abstract}
In this paper, as the second in our series of papers on differential geometry
of microlinear Fr\"{o}licher spaces, we study differenital forms. The
principal result is that the exterior differentiation is uniquely determined
geometrically, just as $\mathrm{div}$ (ergence) and $\mathrm{rot}$ (ation) are
uniquely determined geometrically or physically in classical vector calculus.
This infinitesimal characterization of exterior differentiation has been
completely missing in orthodox differential geometry.

\end{abstract}

\section{Introduction}

Vector analysis is indispensible in studying electromagnetism and fluid
mechanics. The central notions of vector analysis, namely $\mathrm{grad}$,
$\mathrm{div}$ and $\mathrm{rot}$, were introduced infinitesimally as
physically and geometrically meaningful operations. Indeed, their physical or
geometrical meanings determine $\mathrm{grad}$, $\mathrm{div}$ and
$\mathrm{rot}$ uniquely and unequivocally. We know well that the operations
$\mathrm{grad}$, $\mathrm{div}$ and $\mathrm{rot}$ are the precursors of the
exterior differentiation in the modern theory of differential forms.

In a standard course on the theory of differential forms, the exterior
differentiation is introduced by decree as a highly formal operation without
paying little attention to its geometric meaning. Many mathematicians believe
na\"{\i}vely that the proof of Stokes' theorem is easy or even trivial once
the theorem is formulated adequately. We agree completely that its standard
proof is very easy, not to say that it is trivial, but we must insist that the
infinitesimal characterization of the exterior differentiation, which lies at
the core of Stokes' theorem, is not so easy to establish. It is the
infinitesimal Stokes' theorem that underlies the standard (i.e., local or
global) Stokes' theorem. In other words, once the infinitesimal Stokes'
theorem, which is no other than the infinitesimal characterization of the
exterior differentiation, is established, the proof of the standard Stokes'
theorem is highly trivial. This is the whole story of Stokes' theorem from a
conceptual veiwpoint, and its infinitesimal part, which is essential to the
whole story, has been completely missing in orthodox differential geometry.

Although nilpotent infinitesimals are invisible in our standard universe of
mathematics, Weil functors are still meaningful there. The notion of
microlinearity, which is essential to synthetic differential geometry and is
defined in another universe of mathematics, can be externailized by using Weil
functors, as we have discussed in \cite{nishi2}. The principal objective in
this paper is to present the infinitesimal story of Stokes' theorem to
orthodox differential geometers without getting involved seriously in
synthetic differential geometry at all. The reader is strongly recommended to
read Nishimura's \cite{nishi0} as a heuristic introduction to the subject
discussed here.

\section{Preliminaries}

\subsection{Fr\"{o}licher Spaces}

Fr\"{o}licher and his followers have vigorously and consistently developed a
general theory of smooth spaces, often called \textit{Fr\"{o}licher spaces}
for his celebrity, which were intended to be the \textit{underlying set
theory} for infinite-dimensional differential geometry in a sense. A
Fr\"{o}licher space is an underlying set endowed with a class of real-valued
functions on it (simply called \textit{structure} \textit{functions}) and a
class of mappings from the set $\mathbf{R}$ of real numbers to the underlying
set (called \textit{structure} \textit{curves}) subject to the condition that
structure curves and structure functions should compose so as to yield smooth
mappings from $\mathbf{R}$ to itself. It is required that the class of
structure functions and that of structure curves should determine each other
so that each of the two classes is maximal with respect to the other as far as
they abide by the above condition. What is most important among many nice
properties about the category $\mathbf{FS}$ of Fr\"{o}licher spaces and smooth
mappings is that it is cartesian closed, while neither the category of
finite-dimensional smooth manifolds nor that of infinite-dimensional smooth
manifolds modelled after any infinite-dimensional vector spaces such as
Hilbert spaces, Banach spaces, Fr\'{e}chet spaces or the like is so at all.
For a standard reference on Fr\"{o}licher spaces the reader is referred to
\cite{fro}.

\subsection{Weil Algebras and Infinitesimal Objects}

The notion of a \textit{Weil algebra} was introduced by Weil himself in
\cite{wei}. We denote by $\mathbf{W}$ the category of Weil algebras. Roughly
speaking, each Weil algebra corresponds to an infinitesimal object in the
shade. By way of example, the Weil algebra $\mathbf{R}[X]/(X^{2})$ (=the
quotient ring of the polynomial ring $\mathbf{R}[X]$\ of an indeterminate
$X$\ modulo the ideal $(X^{2})$\ generated by $X^{2}$) corresponds to the
infinitesimal object of first-order nilpotent infinitesimals, while the Weil
algebra $\mathbf{R}[X]/(X^{3})$ corresponds to the infinitesimal object of
second-order nilpotent infinitesimals. Although an infinitesimal object is
undoubtedly imaginary in the real world, as has harassed both mathematicians
and philosophers of the 17th and the 18th centuries because mathematicians at
that time preferred to talk infinitesimal objects as if they were real
entities, each Weil algebra yields its corresponding \textit{Weil functor} on
the category of smooth manifolds of some kind to itself, which is no doubt a
real entity. Intuitively speaking, the Weil functor corresponding to a Weil
algebra stands for the exponentiation by the infinitesimal object
corresponding to the Weil algebra at issue. For Weil functors on the category
of finite-dimensional smooth manifolds, the reader is referred to \S 35 of
\cite{kolar}, while the reader can find a readable treatment of Weil functors
on the category of smooth manifolds modelled on convenient vector spaces in
\S 31 of \cite{kri}.

\textit{Synthetic differential geometry }(usually abbreviated to SDG), which
is a kind of differential geometry with a cornucopia of nilpotent
infinitesimals, was forced to invent its models, in which nilpotent
infinitesimals were visible. For a standard textbook on SDG, the reader is
referred to \cite{lav}, while he or she is referred to \cite{kock} for the
model theory of SDG vigorously constructed by Dubuc \cite{dub} and others.
Although we do not get involved in SDG herein, we will exploit locutions in
terms of infinitesimal objects so as to make the paper highly readable. Thus
we prefer to write $W_{D}$\ and $W_{D_{2}}$\ in place of $\mathbf{R}%
[X]/(X^{2})$ and $\mathbf{R}[X]/(X^{3})$ respectively, where $D$ stands for
the infinitesimal object of first-order nilpotent infinitesimals, and $D_{2}%
$\ stands for the infinitesimal object of second-order nilpotent
infinitesimals. To Newton and Leibniz, $D$ standed for
\[
\{d\in\mathbf{R}\mid d^{2}=0\}
\]
while $D_{2}$\ standed for
\[
\{d\in\mathbf{R}\mid d^{3}=0\}
\]
We will write $W_{d\in D_{2}\mapsto d^{2}\in D}$ in place of the homomorphim
of Weil algebras $\mathbf{R}[X]/(X^{2})\rightarrow\mathbf{R}[X]/(X^{3})$
induced by the homomorphism $X\rightarrow X^{2}$ of the polynomial ring
\ $\mathbf{R}[X]$ to itself. Such locutions are justifiable, because the
category $\mathbf{W}$ of Weil algebras in the real world and the category of
infinitesimal objects in the shade are dual to each other. To familiarize
himself or herself with such locutions, the reader is strongly encouraged to
read the first two chapters of \cite{lav}, even if he or she is not interested
in SDG at all.

We need to fix notation and terminology for simplicial objects, which form an
important subclass of infinitesimal objects. \textit{Simplicial objects} are
infinitesimal objects of the form
\begin{align*}
&  D^{n}\{\mathfrak{p}\}\\
&  =\{(d_{1},...,d_{n})\in D^{n}\mid d_{i_{1}}...d_{i_{k}}=0\text{ }%
(\forall(i_{1},...,i_{k})\in\mathfrak{p)\}}%
\end{align*}
where $\mathfrak{p}$\ is a finite set of finite sequences $(i_{1},...,i_{k}%
)$\ of natural numbers between $1$ and $n$, including the endpoints, with
$i_{1}<...<i_{k}$. If $\mathfrak{p}$\ is empty, $D^{n}\{\mathfrak{p}\}$ is
$D^{n}$ itself. If $\mathfrak{p}$ consists of all the binary sequences, then
$D^{n}\{\mathfrak{p}\}$ represents $D(n)$ in the standard terminology of SDG.
Given two simplicial objects $D^{m}\{\mathfrak{p}\}$\ and $D^{n}%
\{\mathfrak{q}\}$, we define a simplicial object $D^{m}\{\mathfrak{p}\}\oplus
D^{n}\{\mathfrak{q}\}$ to be
\[
D^{m+n}\{\mathfrak{p}\oplus\mathfrak{q}\}
\]
where
\begin{align*}
&  \mathfrak{p}\oplus\mathfrak{q}\\
&  =\mathfrak{p}\cup\{(j_{1}+m,...,j_{k}+m)\mid(j_{1},...,j_{k})\in
\mathfrak{q}\}\cup\{(i,j+m)\mid1\leq i\leq m,\ 1\leq j\leq n\}
\end{align*}
Since the operation $\oplus$\ is associative, we can combine any finite number
of simplicial objects by $\oplus$ without bothering about how to insert
parentheses. Given morphisms of simplicial objects $\Phi_{i}:D^{m_{i}%
}\{\mathfrak{p}_{i}\}\rightarrow D^{m}\{\mathfrak{p}\}\ (1\leq i\leq n)$,
there exists a unique morphism of simplicial objects $\Phi:D^{m_{1}%
}\{\mathfrak{p}_{1}\}\oplus...\oplus D^{m_{n}}\{\mathfrak{p}_{n}\}\rightarrow
D^{m}\{\mathfrak{p}\}$ whose restriction to $D^{m_{i}}\{\mathfrak{p}_{i}\}$
coincides with $\Phi_{i}$ for each $i$. We denote this $\Phi$\ by $\Phi
_{1}\oplus...\oplus\Phi_{n}$.

\subsection{Microlinearity}

In \cite{nishi1} we have discussed how to assign, to each pair $(X,W)$\ of a
Fr\"{o}licher space $X$ and a Weil algebra $W$,\ another Fr\"{o}licher space
$X\otimes W$\ called the \textit{Weil prolongation of} $X$ \textit{with
respect to} $W$, which is naturally extended to a bifunctor $\mathbf{FS}%
\times\mathbf{W\rightarrow FS}$, and then to show that the functor
$\cdot\otimes W:\mathbf{FS\rightarrow FS}$ is product-preserving for any Weil
algebra $W$. Weil prolongations are well-known as \textit{Weil functors} for
finite-dimensional and infinite-dimensional smooth manifolds in orthodox
differential geometry, as we have already touched upon in the preceding
subsection. There is a canonical projection from $X\otimes W$ to $X$, and we
denote the inverse image of $x$ under the canonical projection by $\left(
X\otimes W\right)  _{x}$ for any $x\in X$.

The central object of study in SDG is \textit{microlinear} spaces. Although
the notion of a manifold (=a pasting of copies of a certain linear space) is
defined on the local level, the notion of microlinearity is defined absolutely
on the genuinely infinitesimal level. For the historical account of
microlinearity, the reader is referred to \S \S 2.4 of \cite{lav} or Appendix
D of \cite{kock}. To get an adequately restricted cartesian closed subcategory
of Fr\"{o}licher spaces, we have emancipated microlinearity from within a
well-adapted model of SDG to Fr\"{o}licher spaces in the real world in
\cite{nishi2}. Recall that a Fr\"{o}licher space $X$ is called
\textit{microlinear} providing that any finite limit diagram $\mathcal{D}$ in
$\mathbf{W}$ yields a limit diagram $X\otimes\mathcal{D}$ in $\mathbf{FS}$,
where $X\otimes\mathcal{D}$ is obtained from $\mathcal{D}$ by putting
$X\otimes$ to the left of every object and every morphism in $\mathcal{D}$. As
we have discussed there, all convenient vector spaces are microlinear, so that
all $C^{\infty}$-manifolds in the sense of \cite{kri} (cf. Section 27) are
also microlinear.

We have no reason to hold that all Fr\"{o}licher spaces credit Weil
prolongations as exponentiation by infinitesimal objects in the shade.
Therefore we need a notion which distinguishes Fr\"{o}licher spaces that do so
from those that do not. A Fr\"{o}licher space $X$ is called \textit{Weil
exponentiable }if
\begin{equation}
(X\otimes(W_{1}\otimes_{\infty}W_{2}))^{Y}=(X\otimes W_{1})^{Y}\otimes W_{2}
\label{2.1}%
\end{equation}
holds naturally for any Fr\"{o}licher space $Y$ and any Weil algebras $W_{1}$
and $W_{2}$. If $Y=1$, then (\ref{2.1}) degenerates into
\begin{equation}
X\otimes(W_{1}\otimes_{\infty}W_{2})=(X\otimes W_{1})\otimes W_{2} \label{2.2}%
\end{equation}
If $W_{1}=\mathbf{R}$, then (\ref{2.1}) degenerates into
\begin{equation}
(X\otimes W_{2})^{Y}=X^{Y}\otimes W_{2} \label{2.3}%
\end{equation}
We have shown in \cite{nishi1} that all convenient vector spaces are Weil
exponentiable, so that all $C^{\infty}$-manifolds in the sense of \cite{kri}
(cf. Section 27) are Weil exponentiable.

We have demonstrated in \cite{nishi2} that all Fr\"{o}licher spaces that are
microlinear and Weil exponentiable form a cartesian closed category. In the
sequel $M$ is assumed to be such a Fr\"{o}licher space.

\section{Euclidean Vector Spaces}

In this paper we will always mean a preconvenient vector space simply by a
\textit{vector space}. We will choose and fix a vector space $\mathbb{E}$\ in
this sense throughout this section. It is evident that

\begin{lemma}
\label{t3.1}The vector space structure of $\mathbb{E}$\ naturally gives rise
to that of $\mathbb{E}\otimes W$ for any Weil algebra $W$.
\end{lemma}

\begin{proof}
This follows readily from the bifunctionality of $\otimes$\ and the fact that
the functor $\cdot\otimes W:\mathbf{FS}\rightarrow\mathbf{FS}$ is product-preserving.
\end{proof}

\begin{lemma}
\label{t3.2}The vector space structure of $(\mathbb{E}\otimes W_{D}%
)_{\mathbf{0}}$ as the tangent space discussed in our previous paper coincides
with that induced by the vector space structure on $\mathbb{E}\otimes W_{D}$
in the preceding lemma.
\end{lemma}

\begin{proof}
We write $+_{D}$ for the addition in the former vector structure, while we
write $+_{E}$ for the addition in $\mathbb{E}$ as well as for the addition in
$(\mathbb{E}\otimes W_{D})_{\mathbf{0}}$ induced by that in Lemma \ref{t3.1}.
Given $t_{1},t_{2}\in(\mathbb{E}\otimes W_{D})_{\mathbf{0}}$, let us consider
\[
(+_{\mathbb{E}}\otimes\mathrm{id}_{W_{D(2)}})((\mathrm{id}_{\mathbb{E}}\otimes
W_{(d_{1},d_{2})\in D(2)\mapsto d_{1}\in D})(t_{1}),(\mathrm{id}_{\mathbb{E}%
}\otimes W_{(d_{1},d_{2})\in D(2)\mapsto d_{2}\in D})(t_{2}))
\]
It is easy to see that
\begin{align*}
&  (\mathrm{id}_{\mathbb{E}}\otimes W_{d\in D\mapsto(d,0)\in D(2)})\\
&  ((+_{\mathbb{E}}\otimes\mathrm{id}_{W_{D(2)}})((\mathrm{id}_{\mathbb{E}%
}\otimes W_{(d_{1},d_{2})\in D(2)\mapsto d_{1}\in D})(t_{1}),(\mathrm{id}%
_{\mathbb{E}}\otimes W_{(d_{1},d_{2})\in D(2)\mapsto d_{2}\in D})(t_{2})))\\
&  =(+_{\mathbb{E}}\otimes\mathrm{id}_{W_{D(2)}})((\mathrm{id}_{\mathbb{E}%
}\otimes W_{d\in D\mapsto(d,0)\in D(2)})\circ(\mathrm{id}_{\mathbb{E}}\otimes
W_{(d_{1},d_{2})\in D(2)\mapsto d_{1}\in D})(t_{1}),\\
&  (\mathrm{id}_{\mathbb{E}}\otimes W_{d\in D\mapsto(d,0)\in D(2)}%
)\circ(\mathrm{id}_{\mathbb{E}}\otimes W_{(d_{1},d_{2})\in D(2)\mapsto
d_{2}\in D})(t_{2}))\\
&  \text{[by the bifunctionality of }\otimes\text{]}\\
&  =(+_{\mathbb{E}}\otimes\mathrm{id}_{W_{D(2)}})(t_{1},(\mathrm{id}%
_{\mathbb{E}}\otimes W_{d\in D\mapsto0\in D})(t_{2}))\\
&  =t_{1}%
\end{align*}
By the same token, it is also easy to see that
\begin{align*}
&  (\mathrm{id}_{\mathbb{E}}\otimes W_{d\in D\mapsto(0,d)\in D(2)})\\
&  ((+_{\mathbb{E}}\otimes\mathrm{id}_{W_{D(2)}})((\mathrm{id}_{\mathbb{E}%
}\otimes W_{(d_{1},d_{2})\in D(2)\mapsto d_{1}\in D})(t_{1}),(\mathrm{id}%
_{\mathbb{E}}\otimes W_{(d_{1},d_{2})\in D(2)\mapsto d_{2}\in D})(t_{2})))\\
&  =t_{2}%
\end{align*}
Therefore we have
\begin{align*}
&  t_{1}+_{D}t_{2}\\
&  =(\mathrm{id}_{\mathbb{E}}\otimes W_{d\in D\mapsto(d,d)\in D(2)})\\
&  ((+_{\mathbb{E}}\otimes\mathrm{id}_{W_{D(2)}})((\mathrm{id}_{\mathbb{E}%
}\otimes W_{(d_{1},d_{2})\in D(2)\mapsto d_{1}\in D})(t_{1}),(\mathrm{id}%
_{\mathbb{E}}\otimes W_{(d_{1},d_{2})\in D(2)\mapsto d_{2}\in D})(t_{2})))\\
&  =(+_{\mathbb{E}}\otimes\mathrm{id}_{W_{D(2)}})((\mathrm{id}_{\mathbb{E}%
}\otimes W_{d\in D\mapsto(d,d)\in D(2)})\circ(\mathrm{id}_{\mathbb{E}}\otimes
W_{(d_{1},d_{2})\in D(2)\mapsto d_{1}\in D})(t_{1}),\\
&  (\mathrm{id}_{\mathbb{E}}\otimes W_{d\in D\mapsto(d,d)\in D(2)}%
)\circ(\mathrm{id}_{\mathbb{E}}\otimes W_{(d_{1},d_{2})\in D(2)\mapsto
d_{2}\in D})(t_{2}))\\
&  \text{[by the bifunctionality of }\otimes\text{]}\\
&  =t_{1}+_{\mathbb{E}}t_{2}%
\end{align*}

\end{proof}

It is evident that

\begin{proposition}
\label{t3.3}The following conditions on the vector space $\mathbb{E}$\ are equivalent:

\begin{enumerate}
\item The canonical mapping $\mathbf{i}_{\mathbb{E}}:\mathbb{E\times
E\rightarrow E}\otimes W_{D}$ induced by the mapping
\[
(\mathbf{a},\mathbf{b})\in\mathbb{E\times E\mapsto(}x\in\mathbb{R}%
\mapsto\mathbf{a}\mathbb{+}x\mathbf{b}\in\mathbb{E)\in E}^{\mathbb{R}}%
\]
is bijective;

\item The Kock-Lawvere axiom holds in the sense that, for any $t\in
(\mathbb{E}\otimes W_{D})_{\mathbf{0}}$, there exists a unique $\mathbf{a}%
\in\mathbb{E}$ with
\[
t=\mathbf{i}_{\mathbb{E}}(\mathbf{0},\mathbf{a})
\]

\end{enumerate}
\end{proposition}

\begin{definition}
The vector space $\mathbb{E}$\ is called \textit{Euclidean} providing that one
of the above equivalent conditions holds.
\end{definition}

\begin{proposition}
\label{t3.4}If $\mathbb{E}$ is a Euclidena vector space, then so is
$\mathbb{E}^{X}$ for any Fr\"{o}licher space $X$.
\end{proposition}

\begin{proof}
We will check the first condition in Proposition \ref{t3.3}. We have
\begin{align*}
&  \mathbb{E}^{X}\otimes W_{D}\\
&  =(\mathbb{E}\otimes W_{D})^{X}\\
&  =(\mathbb{E\times E})^{X}\\
&  =\mathbb{E}^{X}\mathbb{\times E}^{X}%
\end{align*}
so that we have the desired conclusion.
\end{proof}

\begin{corollary}
The category of Euclidean vector spaces and smooth mappings is cartesian closed.
\end{corollary}

\begin{proposition}
\label{t3.5}If $\mathbb{E}$ is a Euclidena vector space, then so is
$\mathbb{E}\otimes W$ for any Weil algebra $W$.
\end{proposition}

\begin{proof}
We will check the first condition in Proposition \ref{t3.2}. We have
\begin{align*}
&  (\mathbb{E}\otimes W)\otimes W_{D}\\
&  =\mathbb{E}\otimes(W\otimes_{\infty}W_{D})\\
&  =\mathbb{E}\otimes(W_{D}\otimes_{\infty}W)\\
&  =(\mathbb{E}\otimes W_{D})\otimes W\\
&  =(\mathbb{E\times E})\otimes W\\
&  =(\mathbb{E}\otimes W)\times(\mathbb{E}\otimes W)\\
&  \text{[since the functor }\cdot\otimes W\text{ is product-preserving]}%
\end{align*}
so that we have the desired conclusion.
\end{proof}

\begin{remark}
Let $x\in M$. Given $t\in(M\otimes W_{D})_{x}$, we note that $\mathbf{i}%
_{(M\otimes W_{D})_{x}}(\mathbf{0},t)\in(M\otimes W_{D})_{x}\otimes W_{D}$ can
be regarded as an element of $(M\otimes W_{D})\otimes W_{D}=M\otimes W_{D^{2}%
}$, which is no other than
\[
(\mathrm{id}_{M}\otimes W_{(d_{1},d_{2})\in D^{2}\mapsto d_{1}d_{2}\in D})(t)
\]

\end{remark}

We note in passing that

\begin{proposition}
Convenient vector spaces are Euclidean.
\end{proposition}

\begin{proof}
The reader is referred to \S 2 of \cite{kock1}.
\end{proof}

\section{Differential Forms}

Let $\mathbb{E}$ be a Euclidean vector space which is microlinear and Weil exponentiable.

\begin{definition}
Given a smooth mapping $\omega:M\otimes W_{D^{n}}\rightarrow\mathbb{E}$ and a
natural number $i$ with $1\leq i\leq n$, we say that $\omega$ is homogeneous
at the $i$-th position providing that we have
\[
\omega(\alpha\underset{i}{\cdot}\gamma)=\alpha\omega(\gamma)
\]
for any $\gamma\in M\otimes W_{D^{n}}$ and any $\alpha\in\mathbb{R}$, where
$\alpha\underset{i}{\cdot}\gamma$ is defined by
\[
\alpha\underset{i}{\cdot}\gamma=\left(  \mathrm{id}_{M}\otimes W_{\left(
\alpha\underset{i}{\cdot}\right)  _{D^{n}}}\right)  (\gamma)
\]
with the putative mapping $\left(  \alpha\underset{i}{\cdot}\right)  _{D^{n}%
}:D^{n}\rightarrow D^{n}$ being
\[
(d_{1},...,d_{n})\in D^{n}\mapsto(d_{1},...,d_{i-1},\alpha d_{i}%
,d_{i+1},...,d_{n})\in D^{n}%
\]

\end{definition}

\begin{notation}
Given $\gamma\in M\otimes W_{D^{n}}$ and a natural number $i$ with $1\leq
i\leq n$, we denote
\[
\left(  \mathrm{id}_{M}\otimes W_{(d_{1},...,d_{n-1})\in D^{n-1}\mapsto
(d_{1},...,d_{i-1},0,d_{i},...,d_{n-1})\in D^{n}}\right)  \left(
\gamma\right)
\]
by $\gamma\mid_{D^{n-1}}^{i}$.
\end{notation}

\begin{notation}
Given $\eta\in M\otimes W_{D^{n-1}}$ and a natural number $i$ with $1\leq
i\leq n$, we denote
\[
\left\{  \gamma\in M\otimes W_{D^{n}}\mid\gamma\mid_{D^{n-1}}^{i}%
=\eta\right\}
\]
by $\left(  M\otimes W_{D^{n}}\right)  _{\eta}^{i}$.
\end{notation}

\begin{notation}
The putative mapping $\partial_{i}^{n}:D^{n}\rightarrow D^{n}$ is
\[
(d_{1},...,d_{n})\in D^{n}\mapsto(d_{1},...,d_{i-1},d_{n},d_{i},...,d_{n-1}%
)\in D^{n}%
\]

\end{notation}

\begin{remark}
By the natural identification
\[
M\otimes W_{D^{n}}=\left(  M\otimes W_{D^{n-1}}\right)  \otimes W_{D}%
\]
the space $\left(  M\otimes W_{D^{n}}\right)  _{\eta}^{n}$ is a Euclidean
vector space. Under the bijective mapping $\mathrm{id}_{M}\otimes
W_{\partial_{i}^{n}}$, the spaces $\left(  M\otimes W_{D^{n}}\right)  _{\eta
}^{i}$ and $\left(  M\otimes W_{D^{n}}\right)  _{\eta}^{n}$ can be identified,
so that the former is also a Euclidean vector space for any natural number $i
$ with $1\leq i\leq n$. Given $\gamma_{1},\gamma_{2}\in M\otimes W_{D^{n}}$
with
\[
\gamma_{1}\mid_{D^{n-1}}^{i}=\gamma_{2}\mid_{D^{n-1}}^{i}=\eta
\]
we denote the addition of $\gamma_{1}$ and $\gamma_{2}$ in $\left(  M\otimes
W_{D^{n}}\right)  _{\eta}^{i}$ by $\gamma_{1}\underset{i}{+}\gamma_{2}$.
\end{remark}

\begin{proposition}
A smooth mapping $\omega:M\otimes W_{D^{n}}\rightarrow\mathbb{E}$ which is
homogeneous at the $i$-th position is linear at the $i$-th position as well in
the sense that
\[
\omega\left(  \gamma_{1}\underset{i}{+}\gamma_{2}\right)  =\omega(\gamma
_{1})+\omega(\gamma_{2})
\]
for any $\gamma_{1},\gamma_{2}\in M\otimes W_{D^{n}}$ with
\[
\gamma_{1}\mid_{D^{n-1}}^{i}=\gamma_{2}\mid_{D^{n-1}}^{i}%
\]

\end{proposition}

\begin{proof}
The reader is referred to \textit{Proposition 10 in \S 1.2 of \cite{lav}.}
\end{proof}

\begin{definition}
A \textit{differential }$n$\textit{-form }$\omega$\textit{\ on }%
$M$\textit{\ with values in }$\mathbb{E}$ is a smooth mapping
\[
\omega:M\otimes W_{D^{n}}\rightarrow\mathbb{E}%
\]
pursuant to the following conditions:

\begin{enumerate}
\item $\omega$ is $n$-homogeneous in the sense that it is homogeneous at the
$i$-th position for any natural number $i$ with $1\leq i\leq n$.

\item $\omega$ is alternating in the sense that
\[
\omega(\gamma^{\sigma})=\epsilon_{\sigma}\omega(\gamma)
\]
for any permutation $\sigma$ of $1,...,n$, where $\gamma^{\sigma}$ is defined
by
\[
\gamma^{\sigma}=(\mathrm{id}_{M}\otimes W_{\sigma_{D^{n}}})(\gamma)
\]
with the putative mapping $\sigma_{D^{n}}:D^{n}\rightarrow D^{n}$ being
\[
(d_{1},...,d_{n})\in D^{n}\mapsto(d_{\sigma(1)},...,d_{\sigma(n)})\in D^{n}%
\]

\end{enumerate}
\end{definition}

In case that $M$ is a convenient vector space $\mathbb{F}$, we have a more
traditional notion of a differential $n$-form\textit{.}

\begin{definition}
A \textit{differential }$n$\textit{-form}$_{c}$\textit{\ }$\omega$\textit{\ on
}$\mathbb{F}$\textit{\ with values in }$\mathbb{E}$ is a smooth mapping from
$\mathbb{F}$\textit{\ to }$\mathbf{L}_{\mathrm{alt}}^{n}\left(  \mathbb{F}%
;\mathbb{E}\right)  $, where $\mathbf{L}_{\mathrm{alt}}^{n}\left(
\mathbb{F};\mathbb{E}\right)  $ denotes the space of smooth mappings from the
direct product of $n$ copies of $\mathbb{F}$ to $\mathbb{E}$ which are
$n$-linear and alternating.
\end{definition}

\begin{proposition}
In case that $M$ is a convenient vector space $\mathbb{F}$, we assign, to each
\textit{differential }$n$\textit{-form}$_{c}$\textit{\ }$\omega$\textit{\ on
}$\mathbb{F}$\textit{\ with values in }$\mathbb{E}$, the mapping
$\widetilde{\omega}:\mathbb{F}\otimes W_{D^{n}}\rightarrow\mathbb{E}$ with
\[
\widetilde{\omega}\left(  \gamma\right)  =\omega\left(  \mathbf{e}_{1}%
^{\gamma},...,\mathbf{e}_{n}^{\gamma}\right)
\]
for any $\gamma\in\mathbb{F}\otimes W_{D^{n}}$, where
\[
\mathbf{i}_{\mathbb{F}}(\pi(\gamma),\mathbf{e}_{i}^{\gamma})=\left(
\mathrm{id}_{\mathbb{F}}\otimes W_{\iota_{i}^{n}}\right)  \left(
\gamma\right)  \quad\left(  1\leq i\leq n\right)
\]
The assignment gives a mapping from the totality of \textit{differential }%
$n$\textit{-forms}$_{c}$\textit{\ on }$\mathbb{F}$\textit{\ with values in
}$\mathbb{E}$ to that of \textit{differential }$n$\textit{-forms on
}$\mathbb{F}$\textit{\ with values in }$\mathbb{E}$. The mapping is bijective.
\end{proposition}

\begin{proof}
The discussion in Proposition 6 of \S 4.1 in \cite{lav} can be reformulated
easily for our general and abstract context. We should use Proposition 0.3.9
of \cite{nishi-a} in place of Proposition 7 in \S 3.4 of \cite{lav}. The
details can safely be left to the reader, but we note in passing that the
inverse assignment of a \textit{differential }$n$\textit{-form}$_{c}%
$\textit{\ }\underline{$\varrho$}\textit{\ on }$\mathbb{F}$\textit{\ with
values in }$\mathbb{E}$ to each \textit{differential }$n$\textit{-form
}$\varrho$\textit{\ on }$\mathbb{F}$\textit{\ with values in }$\mathbb{E}$
goes as
\[
\underline{\varrho}_{\mathbf{x}}\left(  \mathbf{a}_{1},...,\mathbf{a}%
_{n}\right)  =\varrho\left(  \mathbf{i}_{\mathbb{F}}\left(  \mathbf{x;a}%
_{1},...,\mathbf{a}_{n}\right)  \right)
\]
where $\mathbf{i}_{\mathbb{F}}\left(  \mathbf{x;a}_{1},...,\mathbf{a}%
_{n}\right)  \in\mathbb{F}\otimes W_{D^{n}}$ is the canonical mapping induced
by the mapping
\[
\left(  r_{1},...,r_{n}\right)  \in\mathbb{R}^{n}\mapsto\mathbf{x}%
+r_{1}\mathbf{a}_{1}+...+r_{n}\mathbf{a}_{n}\in\mathbb{F}%
\]

\end{proof}

\section{The Exterior Differentiation}

Let us begin this section with two definitions.

\begin{definition}
Given $\gamma\in\mathbb{E}\otimes W_{D^{n}}$ and a natural number $i$ with
$1\leq i\leq n$, we say that it is homogeneous at the $i$-th position provided
that we have
\[
\left(  \mathrm{id}_{\mathbb{E}}\otimes W_{\left(  \alpha\underset{i}{\cdot
}\right)  _{D^{n}}}\right)  (\gamma)=(\alpha_{\mathbb{E}}\otimes
\mathrm{id}_{W_{D^{n}}})(\gamma)
\]
for any $\alpha\in\mathbb{R}$, where the putative mapping $\left(
\alpha\underset{i}{\cdot}\right)  _{D^{n}}:D^{n}\rightarrow D^{n}$ is
\[
(d_{1},...,d_{n})\in D^{n}\mapsto(d_{1},...,d_{i-1},\alpha d_{i}%
,d_{i+1},...,d_{n})\in D^{n}%
\]
and $\alpha_{\mathbb{E}}$ on the \ right-hand side of the equation stands for
the multiplication by the scalar $\alpha$. We say that $\gamma$ is
$n$\textit{-homogeneous} provided that it is homogeneous at the $i$-th
position for any natural number $i$ with $1\leq i\leq n$.
\end{definition}

\begin{definition}
Given a differential $n$-form $\omega$ on $M$ with values in $\mathbb{E}$ and
$\gamma\in M\otimes W_{D^{n}}$, we define
\[
\int_{\gamma}\omega\in\mathbb{E}\otimes W_{D^{n}}%
\]
as the value of the mapping
\[
\mathbb{E}\times\mathbb{E\cong E}\otimes W_{D}
\begin{array}
[c]{c}%
\mathrm{id}_{\mathbb{E}}\otimes W_{(d_{1},...,d_{n})\in D^{n}\mapsto
d_{1}...d_{n}\in D}\\
\rightarrow\\
\
\end{array}
\mathbb{E}\otimes W_{D^{n}}%
\]
at $(\mathbf{0},\omega(\gamma))\in\mathbb{E}\times\mathbb{E}$.
\end{definition}

It is easy to see that

\begin{proposition}
\label{t5.1}The above mapping
\begin{equation}
\int_{\cdot}\omega:M\otimes W_{D^{n}}\rightarrow\mathbb{E}\otimes W_{D^{n}}
\label{5.1.1}%
\end{equation}
is subject to the following two conditions:

\begin{enumerate}
\item The mapping is a differential $n$-form with values in the vector space
$\mathbb{E}\otimes W_{D^{n}}$;

\item The values of the mapping are all $n$-homogeneous.
\end{enumerate}
\end{proposition}

\begin{proof}
Since the mapping
\[
\mathrm{id}_{\mathbb{E}}\otimes W_{(d_{1},...,d_{n})\in D^{n}\mapsto
d_{1}...d_{n}\in D}:\mathbb{E}\otimes W_{D}\rightarrow\mathbb{E}\otimes
W_{D^{n}}%
\]
preserves the linear structure, we can see readily that the mapping
(\ref{5.1.1}) satisfies the first condition. Since $(\mathbf{0},\omega
(\gamma))$, regarded as an element of $\mathbb{E}\otimes W_{D}$, is
$1$-homogeneous by Lemma \ref{t3.2}, and since the following diagram of
putative mappings
\[%
\begin{array}
[c]{ccccc}
&  & (d_{1},...,d_{n})\in D^{n}\mapsto d_{1}...d_{n}\in D &  & \\
& D^{n} & \rightarrow & D & \\
\left(  \alpha\underset{i}{\cdot}\right)  _{D^{n}} & \uparrow &  & \uparrow &
\alpha_{D}\\
& D^{n} & \rightarrow & D & \\
&  & (d_{1},...,d_{n})\in D^{n}\mapsto d_{1}...d_{n}\in D &  &
\end{array}
\]
commutes, we are sure that the mapping (\ref{5.1.1}) satisfies the second
condition. This completes the proof.
\end{proof}

What is really surprising, we have its converse.

\begin{theorem}
\label{t5.2}If a smooth mapping $\phi:M\otimes W_{D^{n}}\rightarrow
\mathbb{E}\otimes W_{D^{n}}$ abides by the two conditions in Proposition
\ref{t5.1}, then there exists a unique differential $n$-form $\omega$ with
\[
\phi(\gamma)=\int_{\gamma}\omega
\]
for any $\gamma\in M\otimes W_{D^{n}}$.
\end{theorem}

\begin{proof}
The limit diagram of Weil algebras
\[
W_{D}
\begin{array}
[c]{c}%
W_{m_{n}}\\
\rightarrow\\
\qquad
\end{array}
W_{D^{n}}
\begin{array}
[c]{c}%
\ \\
\ \\
W_{\underline{i}_{1}}\\
\rightarrow\\
\vdots\\
\rightarrow\\
\vdots\\
W_{\underline{i}_{n}}\\
\rightarrow\\
W_{0_{n-1}^{n}}\\
\rightarrow
\end{array}
W_{D^{n-1}}%
\]
gives rise to the limit diagram of Fr\"{o}licher spaces
\begin{align*}
&  \mathbb{E}^{M\otimes W_{D^{n}}}\otimes W_{D}
\begin{array}
[c]{c}%
\mathrm{id}_{\mathbb{E}^{M\otimes W_{D^{n}}}}\otimes W_{m_{n}}\\
\rightarrow\\
\qquad
\end{array}
\\
&  \mathbb{E}^{M\otimes W_{D^{n}}}\otimes W_{D^{n}}
\begin{array}
[c]{c}%
\ \\
\ \\
\mathrm{id}_{\mathbb{E}^{M\otimes W_{D^{n}}}}\otimes W_{\underline{i}_{1}}\\
\rightarrow\\
\vdots\\
\rightarrow\\
\vdots\\
\mathrm{id}_{\mathbb{E}^{M\otimes W_{D^{n}}}}\otimes W_{\underline{i}_{n}}\\
\rightarrow\\
\mathrm{id}_{\mathbb{E}^{M\otimes W_{D^{n}}}}\otimes W_{0_{n-1}^{n}}\\
\rightarrow
\end{array}
\mathbb{E}^{M\otimes W_{D^{n}}}\otimes W_{D^{n-1}}%
\end{align*}
because of the microlinearity of $\mathbb{E}^{M\otimes W_{D^{n}}}$, where the
putative mappings $\underline{i}_{j}:D^{n-1}\rightarrow D^{n}$ ($1\leq j\leq
n$) are
\[
(d_{1},...,d_{n-1})\in D^{n-1}\mapsto(d_{1},...,d_{j-1},0,d_{j},...,d_{n-1}%
)\in D^{n}%
\]
while the putative mapping $0_{n-1}^{n}:D^{n-1}\rightarrow D^{n}$ is
\[
(d_{1},...,d_{n-1})\in D^{n-1}\mapsto(0,...,0)\in D^{n}%
\]
Since $\phi$, regarded as an element of $\mathbb{E}^{M\otimes W_{D^{n}}%
}\otimes W_{D^{n}}$ ($=(\mathbb{E}\otimes W_{D^{n}})^{M\otimes W_{D^{n}}}$),
is $n$-homogeneous, it is easy to see that
\[
(\mathrm{id}_{\mathbb{E}^{M\otimes W_{D^{n}}}}\otimes W_{\underline{i}_{1}%
})(\phi)=...=(\mathrm{id}_{\mathbb{E}^{M\otimes W_{D^{n}}}}\otimes
W_{\underline{i}_{n}})(\phi)=(\mathrm{id}_{\mathbb{E}^{M\otimes W_{D^{n}}}%
}\otimes W_{0_{n-1}^{n}})(\phi)
\]
Then the above limit diagram of Fr\"{o}licher spaces guarantees that there
exists a unique $\psi\in\mathbb{E}^{M\otimes W_{D^{n}}}\otimes W_{D}$ with
\[
\phi=(\mathrm{id}_{\mathbb{E}^{M\otimes W_{D^{n}}}}\otimes W_{m_{n}})(\psi)
\]
Since $\mathbb{E}^{M\otimes W_{D^{n}}}$ is Euclidean by Proposition \ref{t3.4}
and $\psi\in(\mathbb{E}^{M\otimes W_{D^{n}}}\otimes W_{D})_{\mathbf{0}}$,
there exists a unique $\omega\in\mathbb{E}^{M\otimes W_{D^{n}}}$ with
\[
\psi=\mathbf{i}_{\mathbb{E}^{M\otimes W_{D^{n}}}}(\mathbf{0},\omega)
\]
Then it is easy to see that $\omega$ is a differential $n$-form with values in
$\mathbb{E}$ such that
\[
\phi(\gamma)=\int_{\gamma}\omega
\]
for any $\gamma\in M\otimes W_{D^{n}}$. This completes the proof.
\end{proof}

\begin{definition}
Given a differential $n$-form $\omega$ on $M$ with values in $\mathbb{E}$, we
define a mapping $\int_{\partial_{i}\cdot}\omega:M\otimes W_{D^{n+1}%
}\rightarrow\mathbb{E}\otimes W_{D^{n+1}}$ to be
\begin{align*}
M\otimes W_{D^{n+1}}
\begin{array}
[c]{c}%
\mathrm{id}_{M}\otimes W_{\partial_{i}^{n+1}}\\
\cong\\
\
\end{array}
M\otimes W_{D^{n+1}}  &  \cong(M\otimes W_{D^{n}})\otimes W_{D}\\%
\begin{array}
[c]{c}%
\left(  \int_{\cdot}\omega\right)  \otimes\mathrm{id}_{W_{D}}\\
\rightarrow\\
\
\end{array}
(\mathbb{E}\otimes W_{D^{n}})\otimes W_{D}  &  \cong\mathbb{E}\otimes
W_{D^{n+1}}%
\end{align*}

\end{definition}

\begin{theorem}
\label{t5.3}Given a differential $n$-form $\omega$ on $M$ with values in
$\mathbb{E}$, there exists a unique differential $(n+1)$-form $\mathbf{d}%
\omega$ on $M$ with values in $\mathbb{E}$ with
\begin{equation}
\int_{\gamma}\mathbf{d}\omega=\sum(-1)^{i+1}\mathbf{D}_{0}\int_{\partial
_{i}\gamma}\omega\label{5.3.1}%
\end{equation}
for any $\gamma\in M\otimes W_{D^{n+1}}$, where $\mathbf{D}_{0}\int
_{\partial_{i}\gamma}\omega$ denotes the mapping
\[
\int_{\partial_{i}\gamma}\omega-(\mathrm{id}_{M\otimes W_{D^{n}}}\otimes
W_{d\in D\mapsto0\in D})\left(  \int_{\partial_{i}\gamma}\omega\right)
\]

\end{theorem}

\begin{proof}
By Theorem \ref{t5.2} it suffices to verify that the right-hand side of
(\ref{5.3.1}) abides by the two conditions in Proposition \ref{t5.1}, which
goes as follows:

\begin{enumerate}
\item We would like to show that $\mathbf{D}_{0}\int_{\partial_{i}\gamma
}\omega\in\mathbb{E}\otimes W_{D^{n+1}}$ is $(n+1)$-homogeneous ($1\leq i\leq
n+1$). Since both\textit{\ }$\int_{\partial_{i}\gamma}\omega$ and
$(\mathrm{id}_{M\otimes W_{D^{n}}}\otimes W_{d\in D\mapsto0\in D}%
)(\int_{\partial_{i}\gamma}\omega)$ are homogeneous at the $j$-th component,
$\mathbf{D}_{0}\int_{\partial_{i}\gamma}\omega$ is homogeneous at the $j$-th
component for $j\neq i$. That $\mathbf{D}_{0}\int_{\partial_{i}\gamma}\omega$
is homogeneous also at the $i$-th component follows from the fact that
$\mathbb{E}\otimes W_{D^{n}}$ is Euclidean.

\item We would like to show that the mapping
\[
\mathbf{D}_{0}\int_{\partial_{i}\cdot}\omega:M\otimes W_{D^{n+1}}%
\rightarrow\mathbb{E}\otimes W_{D^{n+1}}%
\]
is $(n+1)$-homogeneous ($1\leq i\leq n+1$). For $j<i$, it is easy to see that
\begin{align*}
&  \left(  \mathrm{id}_{M}\otimes W_{\partial_{i}^{n+1}}\right)  \circ\left(
\mathrm{id}_{M}\otimes W_{\left(  \alpha\underset{j}{\cdot}\right)  _{D^{n+1}%
}}\right) \\
&  =\left(  \left(  \mathrm{id}_{M}\otimes W_{\left(  \alpha\underset{j}%
{\cdot}\right)  _{D^{n}}}\right)  \otimes\mathrm{id}_{W_{D}}\right)
\circ\left(  \mathrm{id}_{M}\otimes W_{\partial_{i}^{n+1}}\right)
\end{align*}
while, for $j>i$, it is also easy to see that
\begin{align*}
&  \left(  \mathrm{id}_{M}\otimes W_{\partial_{i}^{n+1}}\right)  \circ\left(
\mathrm{id}_{M}\otimes W_{\left(  \alpha\underset{j}{\cdot}\right)  _{D^{n+1}%
}}\right) \\
&  =\left(  \left(  \mathrm{id}_{M}\otimes W_{\left(  \alpha\underset
{j-1}{\cdot}\right)  _{D^{n}}}\right)  \otimes\mathrm{id}_{W_{D}}\right)
\circ\left(  \mathrm{id}_{M}\otimes W_{\partial_{i}^{n+1}}\right)
\end{align*}
Therefore, for $j\neq i$, that
\begin{align*}
&  \left(  \mathbf{D}_{0}\int_{\partial_{i}\cdot}\omega\right)  \circ\left(
\mathrm{id}_{M}\otimes W_{\left(  \alpha\underset{j}{\cdot}\right)  _{D^{n+1}%
}}\right) \\
&  =\left(  \alpha_{\mathbb{E}}\otimes\mathrm{id}_{W_{D^{n+1}}}\right)
\circ\left(  \mathbf{D}_{0}\int_{\partial_{i}\cdot}\omega\right)
\end{align*}
follows directly from the assumption that the mapping
\[
\int_{\cdot}\omega:M\otimes W_{D^{n}}\rightarrow\mathbb{E}\otimes W_{D^{n}}%
\]
is $n$-homogeneous. It remains to show that
\begin{align*}
&  \left(  \mathbf{D}_{0}\int_{\partial_{i}\cdot}\omega\right)  \circ\left(
\mathrm{id}_{M}\otimes W_{\left(  \alpha\underset{i}{\cdot}\right)  _{D^{n+1}%
}}\right) \\
&  =\left(  \alpha_{\mathbb{E}}\otimes\mathrm{id}_{W_{D^{n+1}}}\right)
\circ\left(  \mathbf{D}_{0}\int_{\partial_{i}\cdot}\omega\right)
\end{align*}
which follows readily from
\begin{align*}
&  \left(  \mathrm{id}_{M}\otimes W_{\partial_{i}^{n+1}}\right)  \circ\left(
\mathrm{id}_{M}\otimes W_{\left(  \alpha\underset{i}{\cdot}\right)  _{D^{n+1}%
}}\right) \\
&  =\left(  \mathrm{id}_{M\otimes W_{D^{n}}}\otimes\alpha_{D}\right)
\circ\left(  \mathrm{id}_{M}\otimes W_{\partial_{i}^{n+1}}\right)
\end{align*}
and the Euclideaness of $\mathbb{E}\otimes W_{D^{n}}$.

\item Let $\sigma$ be a permutation of $1,...,n+1$. We would like to show
that
\begin{align}
&  \left(  \sum(-1)^{i+1}\mathbf{D}_{0}\int_{\partial_{i}\gamma}\omega\right)
\circ\left(  \mathrm{id}_{M}\otimes W_{\sigma_{D^{n+1}}}\right)  \nonumber\\
&  =\varepsilon_{\sigma}\sum(-1)^{i+1}\mathbf{D}_{0}\int_{\partial_{i}\gamma
}\omega\label{5.3.2}%
\end{align}
We notice that
\begin{align*}
&  \left(  \mathrm{id}_{M}\otimes W_{\partial_{i}^{n+1}}\right)  \circ\left(
\mathrm{id}_{M}\otimes W_{\sigma_{D^{n+1}}}\right)  \\
&  =\left(  \mathrm{id}_{M}\otimes W_{\left(  \tau_{i}^{\sigma}\right)
_{D^{n}}}\right)  \circ\left(  \mathrm{id}_{M}\otimes W_{\partial_{\sigma
^{-1}(i)}^{n+1}}\right)
\end{align*}
where $\tau_{i}^{\sigma}$ is the permutation of $1,...,n$ with
\begin{align*}
\tau_{i}^{\sigma}(1) &  =\sigma(1),...,\tau_{i}^{\sigma}(\sigma^{-1}%
(i)-1)=\sigma(\sigma^{-1}(i)-1),\\
\tau_{i}^{\sigma}(\sigma^{-1}(i)) &  =\sigma(\sigma^{-1}(i)+1),...,\tau
_{i}^{\sigma}(n)=\sigma(n+1)
\end{align*}
We notice also that
\begin{align*}
&  \left(  \int_{\cdot}\omega\right)  \circ\left(  \mathrm{id}_{M}\otimes
W_{\left(  \tau_{i}^{\sigma}\right)  _{D^{n}}}\right)  \\
&  =\varepsilon_{\tau_{i}^{\sigma}}\int_{\cdot}\omega
\end{align*}
and
\[
\varepsilon_{\tau_{i}^{\sigma}}=(-1)^{\sigma^{-1}(i)-i}\varepsilon_{\sigma}%
\]
Therefore (\ref{5.3.2}) follows.
\end{enumerate}
\end{proof}

\end{document}